\documentclass[a4paper,10pt]{article}
\def \Z {{\mathbf {Z}}}
\def \R {{\mathbf {R}}}
\def \N {{\mathbf {N}}}
\def \T {{\mathbf {T}}}
\def \C {{\mathbf {C}}}
\def\uu{\bigsqcup}
\def\eps{\varepsilon}
\title{   Cпектры cамоподобных эргодических действий}
\author{ Рыжиков В.В.}
\date{14.03. 2022}
\usepackage[T1]{fontenc} % иногда T1 вместо T2A
\usepackage[cp1251]{inputenc}  
\usepackage[russian]{babel}
\usepackage{graphicx}
\graphicspath{{pictures/}}
\DeclareGraphicsExtensions{.pdf,.png,.jpg}

\begin{document}
%\udk{517.9}

\maketitle

\begin{abstract} 
We consider self-similar constructions of transformations that preserve the sigma-finite measure.
  Their properties and spectra of induced Gaussian and Poisson dynamical systems are studied.
An orthogonal operator corresponding to such a transformation has the following property: 
some of its powers is a nontrivial direct sum of operators isomorphic to the original one.
 For any set $M\subset \N$, there is a  Poisson suspension with   spectral multiplicity set of the form $M\cup\{\infty\}$.
There is  Gaussian flow $S_t$ such that  the automorphism $S_{p^{n}}$ has  spectral multiplicities $\{1,\infty\}$ as $n\leq 0$, and $\{p^n,\infty\}$ as $n> 0$, its infinite tensor power $ T_t$  is so that
automorphisms $T_{p^{n}}$ have different spectral types for $n\leq 0$, but
 all automorphisms $T_{p^{n}}, n>0,$ are  isomorphic.  We  also consider the constructions of self-similar flows and discuss the problems that arise in connection with them.

\vspace{3mm}
В заметке рассматриваются самоподобные конструкции  преобразований, сохраняющих сигма-конечную меру,
  изучаются их свойства и спектры индуцированных  гауссовских и пуассоновских динамических систем.
Ортогональный  оператор, отвечающий такому преобразованию, обладает свойством:  некоторая его степень   является нетривиальной  прямой суммой  операторов, изоморфных исходному. 
Получены следующие результаты. 
 Для любого подмножества $M$ натурального ряда в классе пуасоновских надстроек реализованы  наборы спектральных кратностей вида $M\cup\{\infty\}$.  
Предъявлен  гауссовский поток $S_t$  такой, что  автоморфизмы $S_{p^{n}}$ обладают набором спектральных кратностей $\{1,\infty\}$, если   $n\leq 0$,  и наборами кратностей 
$\{p^n\infty\}$ при $n> 0$;     бесконечная  тензорная степень $ T_t$ потока $S_t$ такова, что 
автоморфизмы $T_{p^{n}}$ обладают различными спектральными типами при   $n\leq 0$, но 
 все автоморфизмы $T_{p^{n}}, n>0,$ попарно изоморфны между собой.

\end{abstract}
%\begin{keywords}
  %  Сохраняющие меру преобразования, самоподобные конструкции,  слабое замыкание,  спектр, изоморфизм эргодических систем.
%\end{keywords}
\large
\section{Введение}
Преобразования ранга 1 находят разнообразные приложения
 в эргодической теории. Например,  связи с известными в  спектральной теории динамических систем  проблемами Колмогорова и Рохлина   были указаны  явные конструкции перемешивающих преобразований ранга 1, тензорный квадрат которых обладает  однородным непростым спектром, а спектральные меры  этих преобразований не обладают групповым свойством  \cite{R20}.  
Для  ответа на  с вопрос М. И. Гордина о существовании преобразования  с простным сингулярным спектром и эргодической гомоклинической группой  также привлекались  преобразования ранга 1 \cite{R14}. 
Действия ранга 1 применялись  в  решении  задачи В. Бергельсона о сосуществовании жестких и перемешивающих последовательностей для сохраняющего меру $Z$-действия \cite{R21}   и  задачи В. И. Оселедца о сингулярных и абсолютно непрерывных распределениях случайной величины вида $\xi(x)+a\xi(y)$ в зависимости от параметра $a$ (см. \cite{R22}).  
 Настоящая  заметка посвящена  самоподобным конструкциям ранга 1 и их приложениям к  спектральной теории динамических систем.    Ортогональные операторы, отвечающие  самоподобным действиям,  индуцируют    необычные  спектральные свойства  гауссовских действий  и пуассоновских надстроек.

Напомним, что гауссовские системы можно рассматривать как результат   вложения ортогональной группы $O(\infty)$ в группу автоморфизмов
  пространства $ R^\infty$ с гауссовской мерой.
 Пуассоновские  системы  появляются благодаря действию   группы автоморфизмов пространства с сигма-конечной мерой    на конфигурационном пространстве с вероятностной мерой Пуассона  (см., например,\cite{KSF},\cite{N}). Таким образом, с преобразованием   $ T$ пространства  с сигма-конечной мерой ассоциированы два замечательных объекта: гауссовский автоморфизм $ G(T)$ и пуассоновская надстройка  $ P(T)$. Нам в дальнейшем понадобится  следующий факт: операторы 
$ G(T)$ и  $ P(T)$ изоморфны  оператору
$ \exp(T) =\bigoplus_{n = 0}^{\infty} T^{\odot n},$
где $ T^{\odot 0}$ -- одномерный тождественный оператор,
$ T^{\odot n}$ - симметрическая тензорная $n$-степень оператора $T$. 

Предположим, что унитарные операторы $T_i$, $i \in\N$   имеют простой спектр и   попарно  дизъюнктны, 
 т.е. их спектральные меры взаимно сингулярны. Обозначим через  ${\bf m} T$  прямую сумму $m$ копий оператора $T$. 
Если всевозможные тензорные произведения операторов
$T_{i_1}\otimes\dots\otimes T_{i_k}$, $i_m\in \N$  дизъюнктны с операторами $T_i$, $i\in \N$,  и имеют счетнократный спектр,
то  $M\cup\{\infty\}$  является множеством спектральных кратностей оператора 
$$   \exp\left( \bigoplus_{m\in M} {\bf m} T_{m}   \right).$$

Мы предъявим семейство  $T_i$, $i\in \N$ самоподобных  преобразований, которые как операторы будут удовлетворять  описанным выше свойствам (преобразования и индуцировнные ими операторы в статье обозначаются одинаково). Это приводит к следующему результату.

\vspace{2mm}
\bf Теорема 1. \it Для любого подмножества $M$ натурального ряда найдется  пуасcоновская надстройка, обладающая  набором спектральных кратностей  $M\cup\{\infty\}$. \rm

\vspace{2mm}
Неизвестно, всегда ли найдется  эргодическое преобразование вероятностного пространства с  набором спектральных кратностей  вида $M\cup\{m\}$  при $m>2$.
В классе эргодических преобразований реализованы  всевозможные  наборы кратностей вида $M\cup\{1\}$,  $M\cup\{2\}$  (см. \cite{A} и обзоры \cite{D}, \cite{KL}).

Унитарный оператор   c простым спектром (иначе говоря, с кратностью  спектра, равной   1)
  эквивалентен оператору  умножения на переменную $z$, 
$|z|=1$, действующему  
в пространстве $L_2(\T,\sigma)$. Здесь $\sigma$ -- борелевская мера на единичной 
окружности $\T=\{ t:  |t|=1\}$, которую называют спектральной мерой оператора.   Если  $\sigma$ симметрична относительно инволюции $z\to\bar z$,  то соответствующий оператор умножения   является также ортогональным.

Пусть унитарный оператор $T$ с простым спектром удовлетворяет свойству:
степень $T^p$ подобна  оператору  $\bigoplus_{i=1}^p U$ для некоторого   $p>1$.
Назовем такой оператор $p$-самоподобным, а его спектральную меру назовем $p$-мерой.
Несложно заметить, что  $p$-мера квазиинвариантна относительно
стандартного действия полициклической группы 
$$C_p=\left\{\exp\left(\frac {2\pi i k} {p^n}\right) \: \ k\in \Z,  n\in \N\right\}$$ 
на группе $\T$.   Располагая части  $p$-мер на прямой $\R$, 
можно конструировать  потоки в группе $O(\infty)$, которые  индуцируют   гауссовские  потоки. 
Используя этот подхол,  мы докажем следующее утверждение.

\vspace{2mm}
\bf Теорема 2. \it Найдется   гауссовский поток $S_t$ такой, что автоморфизмы $S_{p^{n}}$ 
обладают набором спектральных кратностей $\{1,\infty\}$ для всех    $n\leq 0$,  и наборами кратностей 
$\{p^n,\infty\}$ при $n> 0$.     Для некоторого  гауссовского потока $T_t$
автоморфизмы $T_{p^{n}}, \ n\leq 0$ обладают различными спектральными типами, но 
  автоморфизмы $T_{p^{n}},\ n>0,$ попарно изоморфны между собой.

\rm

\vspace{2mm} 
\bf Замечание. \rm Отметим, что С.В.Тихонов, используя результаты \cite{S},  привел пример негауссовского потока  $T_t$, для которого $T_{\frac 1 2}$ не изоморфнен    $T_{1}$, но $T_{ 1}$ изоморфен $T_{2}$ (см. \cite{DR}).  Типичный автоморфизм вероятностного пространства вкладывается в континуум неизоморфных потоков \cite{S} и этим  он похож на  гауссовские автоморфизмы с простым спектром, 
так как последние  также включаются  в континуальное множество спектрально неизоморфных потоков.

\section{Самоподобные конструкции  и их свойства}
Конструкция преобразования ранга 1 (см., например, в \cite{R20}) определяется параметрами 
 $h_1$, $r_j=r$ и 
$$ \bar s_j=(s_j(1), s_j(2),\dots,s_j(r_j)), \ s_j(i)\geq 0, \ r_j>1. $$
Положив  $h_1=h$, $r_j=r$,  
$$ \bar s_j=(s(1)h_j, s(2)h_j,\dots,s(r)h_j), \ \ h_{j+1}= q^j h, $$ 
 получим так называемую самоподобную конструкцию с коэффициентом подобия $q=r+s(1)+ s(2)+\dots +s(r)$. 

Для удобства читателя мы даем полное описание конструкции в частном случае для параметров
$$h_1 = h, \ r = 2, \ s(1) = 1, \ s(2) = p-3.  $$
  Унитарный оператор, отвечающий такой конструкции,  будет $p$-самоподобным. 
Выбор параметров $s(1)=1, s(2)=p-3 $ не играет принципиальной роли, он связан с тем, 
что в \cite{KR} была рассмотрена такая конструкция $T$
при $p=8$, $h=1$ и описаны все ее слабые пределы. Мы этим воспользуемся позже.

\bf Конструкция типа (h, p). \rm  Преобразование $T$ и его фазовое пространство $X$ строятся индуктивно.  
Пусть на  шаге $j\geq 1$ определен
набор из $h_j$ непересекающихся интервалов одинаковой длины:
$$E_j, TE_j, T^2E_j,\dots, T^{ h_j-1}E_j,$$
на интервалах  $E_j,  \dots, T^{ h_j-2}E_j$
преобразование $T$ является  обычным   переносом.  

Такой набор  называют башней высоты $h_j$. 
Основание башни $E_j$ разбиваем два интервала  $E_j^1,E_j^2$ одинаковой меры. 
  Башню разбиваем  на две колонны: первая состоит из интервалов 
$$E_j^1, TE_j^1, T^2E_j^1,\dots, T^{ h_j-1}E_j^1,$$
а вторая -- из интервалов 
$$E_j^2, TE_j^2, T^2E_j^2,\dots, T^{ h_j-1}E_j^2.$$
Добавляем $h_j$ новых интервалов
$$T^{h_j}E_{j}^1, T^{h_j+1}E_{j}^1,\dots, T^{2h_{j}-1}E_{j}^1$$
 над первой колонной и добавляем $(p-3)h_j$  интервалов
$$T^{3h_j}E_{j}^2, T^{h_j+1}E_{j}^2,\dots, T^{ph_{j}-1}E_{j}^2$$
 над второй колонной. Отметим, что по построению все добавленные интервалы не пересекаются со старыми.
  Положив 
$$T^{ 2h_j}E_j^1 =E_j^2, \ E_{j+1}= E_j^1,$$
получаем башню этапа $j+1$, состоящую  из $h_{j+1}=ph_j$ этажей
$$E_{j+1}, TE_{j+1}, T^2 E_{j+1},\dots, T^{h_{j+1}-1}E_{j+1}.$$
Преобразование $T$    действует на этой башне как  обычный перенос  интервалов, пока оно  не определено на последнем интервале. 
Доопределяя преобразование, мы  не меняем того, что было задано на предыдущих этапах.
  Продолжая построение до бесконечности,  получим обратимое преобразование $T$  
  множества 
 $$   X =\bigcup_j \bigcup_{i=0}^{h_{j}-1} T^{i}E_{j},$$
сохраняющее меру Лебега.
Известно,  что такое преобразование эргодично и   имеет простой спектр (как и все преобразования ранга 1).

\bf Самоподобие. \rm 
Замечаем, что  множество  
$$  \tilde X =\bigcup_j \bigcup_{i=0}^{h_{j}-1} T^{pi}E_{j+1}$$
инвариантно относительно степени $T^p$.  
  Мера $\mu$ имеет $p$ эргодических компонент относительно преобразования $T^p$, 
эти компоненты сосредоточены на непересекающихся инвариантных множествах $\tilde X, T\tilde X,\dots, T^{p-1}\tilde X$.
Ограничение степени  $T^p$ на любое из   инвариантных множеств $T^k\tilde X$ будет  подобно преобразованию $T$. 
Подобие между $T^p|\tilde X$ и $T$  осуществляется отображением, которое увеличивает меру в $p$ раз, сопоставляя интервалу 
$T^{pi}E_{j+1}$  интервал $T^iE_j$ при $0\leq i<h_j$.  Можно сказать, что на $\tilde X$ преобразование 
$T^p$ устроено точно так, как $T$ на $X$,  только теперь роль интервалов $T^iE_j$ играют интервалы  $T^{pi}E_{j+1}$.
Оператор $T$ обладает   $p$-самоподобием  и   набором  других свойств, которые   
рассмотрены  ниже. 

\vspace{2mm}
\bf Теорема 3. \it  Пусть $T$ --  преобразованиe типа $(m,p)$, а числа 
$q\neq q'$  взаимно просты с $p$.  Тогда  степени $T^q$ и $T^{q'}$ имеют простые взаимно сингулярные спектры.\rm 

 \vspace{2mm}
\bf Теорема 4. \it  Пусть $S,T$  являются  преобразованиями типов $(m,p)$ и $(n,p)$, соответственно,
причем  $m,n,p$   взаимно просты.  Тогда $S,T$  при $p\geq 8$ обладают  следующими  свойствами. 

 {(1)}  Ненулевые слабые пределы степеней оператора $S$ (аналогично для $T$) имеют вид $\frac 1 {2^m} S^q$, где  $m\in \N$ 2, $q\in\Z$. 

 (2) Cпектральная мера преобразования     $S$ сингулярна  
 и квазиинвариантна относительно действия  группы  $C_p$   на    $\T$. 

(3)  $S$ и $T$  спектрально дизъюнктны.

(4) Спектры произведений  $T\otimes T$ и $T\otimes S$   счетнократны.
\rm

\section{Доказательства теорем}
Доказательство теоремы 4.  (1). Для конструкции типа (1,8)  утверждение (1)   установлено в \cite{KR}. 
Общий случай доказывается совершенно аналогично.

\vspace{2mm}
(2).  Слабые пределы $S^{-h_j}\to \frac 1 2 I$ обеспечивают сингулярность спектральной меры $\sigma$.  Спектр $S$ простой, а спектр $S^p$ однородный кратности $p$. Это означает, что на дугах $\left[\frac {2\pi (k-1)} p,\frac {2\pi k} p\right)$, $k=1,\dots,p$, 
ограничение меры $\sigma$ при отображении $z^p$ переходит в меру, эквивалентную мере $\sigma$. Тогда ограничение 
меры $\sigma$ на дуги 
$\left[\frac {2\pi (k-1)} {p^2},\frac {2\pi k} {p^2}\right)$ также подобны мере $\sigma$.  Получаем окончательно, что 
отображение $z^{p^n}$ осуществляет подобие всех ограничений меры $\sigma$ на дугах 
$\left[\frac {2\pi (k-1)} {p^n},\frac {2\pi k} {p^n}\right)$, $k\in \N$,  с мерой $\sigma$ на окружности $\T$.
Таким образом, поворот окружности $\T$ на угол $\frac {2\pi k} {p^n}$, $k\in \N$, переводит меру $\sigma$  в эквивалентную ей. Это означает квазиинвариантность $\sigma$  относительно действий полиициклической группы $C_p$.

\vspace{2mm}
(3). Пусть $SJ=JT$, тогда $S^{n}J=JT^{n}$.  Если  $S^{n_k}\to \frac 1 2 I$, $n_k>0$, то, как показано в $\cite{KR}$, для всех больших $k$ выполнено
$n_k=h_{j_k}=mp^{j_k}$. Таким образом, имеем следующие слабо сходящиеся последовательности:
$$ S^{n_k}J\to \frac 1 2 J, \ \ JT^{n_k}\to \frac 1 2 J,$$
$$ S^{mp^{j_k}}J\to \frac 1 2 J, \ \ JT^{mp^{j_k}}\to \frac 1 2 J.$$
Операторы вида $ \frac 1 2 T^s$ не могут быть предельными точками для множества операторов $T^{mp^{j_k}}$,
так как они являются предельными точками только для множеств содержащих операторы  $T^{np^{j_k}+s}$, но множества 
$\{mp^{j}:j\in\N\}$ и  $\{np^{i}+s:i\in\N\}$ имеют лишь конечные пересечения
в силу того, что $m,n,p$ взаимно просты. Таким образом, слабыми предельными точками операторов $T^{mp^{j_k}}$
могут быть только 0 или операторы вида $ \frac 1 {2^n} T^s$ при $n>1$. Но тогда из $JT^{mp^{j_k}}\to \frac 1 2 J$
очевидным образом вытекает $J=0$.

\vspace{2mm}
(4). Мера $\sigma_S\times\sigma_T$ (и мера $\sigma_S\times\sigma_S$) на двумерном 
торе квазиинвариантна относительно группы сдвигов на векторы вида
$\left(\frac {2\pi k} {p^n}, \frac {-2\pi k} {p^n}\right)$. Это означает, что мера $\sigma\times\sigma$
проектируется на диагональ $\{(z,z):  z\in \T\}$ в торе с бесконечной кратностью, так как  
соответствующие условные меры на слоях, ортогональных диагонали, не могут быть конечной суммой точечных мер.
Следовательно, спектр тензорных произведений  $p$-самоподобных операторов является счетнократным. Теорема 4 доказана.

\bf Доказательство теоремы 3. \rm
Рассмотрим частный случай,  когда  $T$  -- преобразованиe типа $(1,p)$.
Известно, что $T^1$ как преобразование ранга 1 имеет простой спектр.
Покажем, что при $q\perp p$ (нет общего делителя)  степень $T^q$ имеет простой спектр.
Найдется последовательность  $qn_i=p^{i} +s,$  $n_i\to\infty$,  такая, что  $s<q$.  
Если $s=p^ks'$, разделим выражения в  равенстве $qn_i=p^{i} +s$  на  
 $p^k$. Таким образом,  считаем, что $s$ не делится на $p$.
Так как 
$$  (T^q)^{n_i}\to \frac 1 2 T^s,$$
простота спектра  $T^s$ влечет за собой простоту спектра степени $T^q$ (циклический вектор для $T^s$ является 
циклическим для $T^q$).  В предположении, что $T^q$ имеет непростой спектр, мы показали
отсутствие  минимального  $s\perp p$ такого, что $T^s$ имеет непростой спектр. Следовательно, предположение неверно,
а $T^q$ имеет простой спектр.

Покажем, что $T^q$ и  $T^r$ при $q\perp r$ не имеют ненулевых сплетений.   Предположим, что  $T^qJ=JT^r$ и $J\neq 0$.  Рассмотрим последовательность вида  $qn_i=p^{i} -k_i$, где $|k_i|<q$.  Далее рассматриваем, только такие $i$, когда  
$k_i=k$. Имеем
$$ T^{qn_i} J=JT^{rn_i},\ \  \frac {T^kJ} 2 =JP,\ \ \frac{\|J\|} 2 =\|JP\|, $$
где $P$ -- слабая предельная точка для  степеней $T^{rn_i}$.
Условие $\|JP\|= \frac{\|J\|} 2$ означает, что $P=\frac {T^mJ} 2$ (другие слабые пределы не подходят).
Тогда, сравнивая  $rn_i=p^{j_i} - m$  с $qn_i=p^{i} -k$, получаем 
$q=p^{n}r$ для некоторого целого $n\neq 0$. Так как $q\perp p$ и $q\perp r$, получаем противоречие.
Следовательно, что $J=0$.

\bf Доказательство теоремы 1. \rm
Пусть  $p_m$, $m\in M$, являются  различными  нечетными простыми числами, 
через  $T_m$ обозначим преобразование типа $(p_m, 8)$.  Покажем, что     $M\cup\{\infty\}$   
есть множество спектральных кратностей оператора 
$$   \exp\left( \bigoplus_{m\in M} {\bf m} T_{m}   \right)=
1\oplus \left(\bigoplus_{m\in M} {\bf m} T_{m}\right)\oplus  Q.$$
Замечаем, что  
$$Q= \left(\bigoplus_{m\in M} {\bf m} T_{m}\right)^{\odot 2}\oplus
 \left(\bigoplus_{m\in M} {\bf m} T_{m}\right)^{\odot 3}\oplus 
\dots$$
имеет счетнократный спектр в силу пункта (4) теоремы 3.
Спектры операторов $T_m$  взаимно сингулярны со спектром  оператора $Q$. Действительно, 
степени операторов $T_m$ имеют предельную точку вида $\frac 1 2 I$, причем нормы
всех слабых пределов не превосходят $\frac {1} 2$.
   Оператор $Q$ не может иметь нетривиальное спелетение с $T_m$, покажем это.
Пусть $T_mJ=JQ$. Тогда  $T_m^nJ=JQ^n$, переходя к слабым пределам, кода $T_m^{n_k}\to \frac I 2$, 
 получим
$\frac J 2 = J R,$ где $R$ -- слабый предел степеней оператора $Q$.
Но нормы слабых пределов  $Q^n$ не превосходят $\frac {1} 4$. Имеем $\frac {\|J\|} 2 \leq \frac {\|J\|} 4$,
следовательно,  $J=0$.
Таким образом, $T_m$ не имеют нетривиальных сплетений с оператором $Q$.

Набор спектральных кратностей  оператора $P_1=\bigoplus_{m\in M} {\bf m} T_{m}$  есть $M$, что вытекает из (3) теоремы 3
(кратность 1, которую дает тождественный оператор в пространстве констант,  по традиции не учитывается).
Спектральный тип оператора $P_1$ дизъюнктен со спектром оператора $Q$. Получили, что $M\cup\{\infty\}$   является 
набором спектнальных кратностей оператора $\exp(P_1)$, что и требовалось. Теорема 1 доказана.

\bf Доказательство теоремы 2. \rm  Рассмотрим ортогональный поток $O_t$ на $L_2(\R,\sigma_1)$, где мера $\sigma_1$ 
совпадает
с мерой $\sigma$ при  отождествлении  окружности $\T$  с полуинтервалом  $\left[-\frac \pi 2, \frac \pi 2\right)$  
посредством функции  $arg(z)$. Гауссовский поток $S_t=\exp(O_t)$ обладает искомыми свойствами. 
Это вытекает из того, что $O_{p^n}$ имеет простой спектр при $n\leq 0$ и однородный кратности $p^n$ при $n> 0$.
Симметрические тензорные степени оператора $O_t$ имеют счетнократный спектры,  взаимно сингулярные со спектром $O_t$.

Если в качестве $T_t$ взять бесконечную тензорную степень потока $\exp(O_t)$, 
автоморфизмы $T_t$ $0<t\leq 1$ обладают различными спектральными типами, но 
  автоморфизмы $T_{p^{n}}$ при  $ n\geq 0,$ имеют одинаковый счетнократный спектр и попарно изоморфны между собой,
что легко проверяется.   Теорема 2 доказана.

\section{Самоподобные кострукции потоков}
 Поток  ранга 1 задается параметрами 
 $h_1,w> 0$,  $h_1,w\in\R$, натуральными $r_j>1$ и векторами 
$$ \bar s_j=(s_j(1), s_j(2),\dots,s_j(r_j)),  \ s_j(i)\in \R^+. $$
В построении фазового пространства потока на этапе $j$ башня отождествляются с 
прямоугольниками высоты $h_j$ и ширины $w$.
Башня разрезается на $r_j$ колонн (одинаковые прямоугольники высоты $h_j$), над колоннами надстраиваются прямоугольники высотой $s_j(i)$ (их ширина равна ширине колонн),
затем верх надстроенной $i$-той колонны склеивается с низом $i+1$ колонны. В результате мы получаем
башню этапа $j+1$, которую отождествляем теперь с прямоугольником высоты $h_{j+1}=rh_j +\sum_{i=1}^r s_j(i)$ 
ширины $w/r_1r_2\dots r_j$.   Движение точки под действием потока в башне есть движение  вверх с постоянной скоростью. При достижении верхней границы башни этапа $j$,
движение точки рассматривается  в башне этапа $j+1$ и так далее. Объединение всех башен является фазовым пространством
построенного потока.

Рассмотрим простейшую $q$-самоподобную конструкцию потока $T_t$ с параметрами   $h_1=1$, $w=1$, $r_j=2$,  
$ \bar s_j=(0, (q-2)h_j)$, где  $q>2$,  $h_{j+1}= q^j h$.
Подобие потоков $T_t$ и $T_{qt}$  устанавливается  следующим образом.

Фазовое пространство $X$ является объединением башен $X_j$, причем
башня $X_{j+1}$ в $q$ раз выше и в 2 раза
уже башни $X_j$. Обозначим через $\Phi_j:X_{j+1}\to X_j$ соответствующее отображение  подобия прямоугольников $X_{j+1}$ и $X_j$ (оно не сохраняет меру, так как $q>2$. 
Теперь определим самоподобие пространства
$\Phi_j:X\to X$, положив $\Phi|_{X_{j+1}}=\Phi_j|_{X_{j+1}}$. 
Чтобы явно увидеть самоподобие пространства и потока, удобно 
фазовое пространство  $X$  представить как части плоскости, изображенные слева и справа на 
рисунке. Правая часть ассоциируется с объединением $\cup_{j=1}^\infty X_j$,  а левая  
 с  $\cup_{j=1}^\infty \hat X_{j}$, где  $\hat X_{j}$ -- перерисованная башня $X_{j+1}$.
 
\vspace{3mm}
\begin{center}
\includegraphics{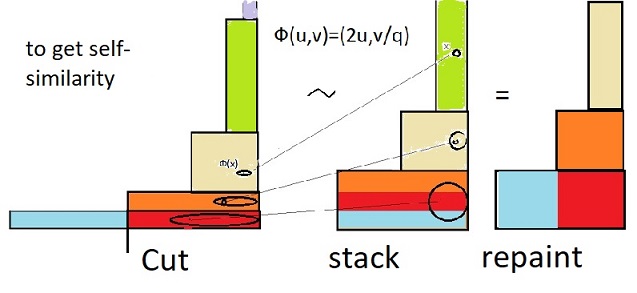}
\end{center}

Подобие $\Phi$, заданное отображением     $\Phi(u,v)=(2u,v/q)$, сопрягает поток $T_t$ с потоком
$T_{qt}$. Гауссовская надстройка над таким потоком пополняет коллекцию  потоков,
самоподобных в смысле работы \cite{DR}. Пуассоновская надстройка наследует спектральное самоподобие.

\vspace{3mm}
\bf Замечания и вопросы. \rm Обратим внимание читателя на следующее утверждение.
\it  Если  унитарный поток $U_t$ обладает простым сингулярным спектром, то для почти всех $t$  операторы  $U_t$  имеют  простой спектр.  \rm  Наш поток $T_t$ имеет  сингулярный 
спектр, так как сверточные степени спектральной меры потока взаимно сингулярны. Последнее
вытекает из наличия слабого предела $T_{q^n}\to_w \frac I 2$.
  Так как спектр потока (ранга 1)  простой,   возникает 
вопрос: является ли  простым  спектр автоморфизма $T_1$в случае  иррационального $q$? 
Если  $q$ -- целое число,  спектр   $T_1$ сингулярный и счетнократный. 
(Пуассоновская надстройка $P(T_t)$ для рациональных $t$  имеет счетнократный спектр,
но для почти всех $t$  в набор спектральных кратностей содержит 1.) 
Если   дробные части последовательности $q^n$ имеют предельную иррациональную точку $\alpha$,
то  несложно показать, что $T_1$ имеет простой спектр, применяя возникающие в этом случае  нетривиальные слабые пределы
вида $$\frac {T_{\{k\alpha\}}} {2^k}, \ \ k=1,2,\dots.$$
Интересная  нерешенная задача о спектральной кратности $T_1$  возникает в случае, когда $q$ является числом Пизо (Pisot number), т.е. расстояние от $q^n$ до множества $\N$ стремится к 0. 

Задача  исследования спектрального типа и спектральной кратности 
тензорных произведений вида $T_t\otimes T_{at}$ в зависимости от $q\in (2,+\infty)$ и
$a\in [1,+\infty)$ также представляет определенный интерес, особенно для небольших значений $q$.

\newpage
%\section

\end{document}